\pdfoutput=1

\documentclass[11pt,a4paper%
,fleqn%
,DIV=13%
,BCOR=10mm%
,twoside%
,headings=normal%
,numbers=enddot%
,headsepline=on%
]{scrartcl}

\input zbm_def.input

\title{The longest perimeter small triacontadigon}
\author{Bernd Mulansky\footnote{Institute of Mathematics,
Clausthal University of Technology, Erzstr.~1, 38678 Clausthal-Zellerfeld, Germany
(\texttt{bernd.mulansky@tu-clausthal.de})}
\and Andreas Potschka\footnote{Institute of Mathematics,
Clausthal University of Technology, Erzstr.~1, 38678 Clausthal-Zellerfeld, Germany
(\texttt{andreas.potschka@tu-clausthal.de})}
}
    
\date{}

\begin{document}

\maketitle

\abstract{We showcase a small triacontadigon, \ie a convex $32$-gon of unit diameter,
with perimeter $3.1403311569546$, which improves on all examples available so far.}

\section{Introduction}
In 1922, \Name{Karl Reinhardt} \cite{Reinhardt1922} published his results
on isodiametric extremal problems for convex polygons. He asked for
the maximal area or the maximal perimeter of convex $n$-gons of
unit diameter, which nowadays are usually called \notion{small polygons}.
The following theorem summarizes his findings on the longest perimeter. 
\begin{theorem}[Reinhardt 1922] \label{th:lpsp}
	Let $n\ge 3$ be an integer. The perimeter of a  small convex $n$-gon 
	is always less or equal to $2n \sin\frac\pi{2n}$.
	Equality is attainable if and only if $n$ is not a power of $2$. 
	For odd $n$ the regular small $n$-gon has longest perimeter,
	but for even $n$ the regular small $n$-gon is never extremal.
\end{theorem}
\Name{Reinhardt} also gave a fairly complete description of the 
extrema in case that $n$ has an odd factor, 
thereby noticing firstly that the collections of extreme 
$n$-gons depend on the divisibility properties of $n$.

For surveys on several variations of isodiametric extremal problems
for polygons we refer to \cite{Audet2007,Audet2009a,Mossinghoff2006a,Toth2022}.

Here we focus on the maximal perimeter of a small $n$-gon when $n$ is a power
of $2$. The solution for $n=4$ is usually attributed to \cite{Tamvakis1987,Datta1997}. 
However, the extreme quadrilateral was already derived by \Name{Taylor}
\cite[Corollary on p.~193]{Taylor1953},
referenced in \cite[2.25 on p.~456]{Mitrinovic1989}.
A detailed proof is also contained in
\cite[86c on p.~279--284]{Shklyarskii1974}, which we could
find due to a hint in \cite{Gashkov1985}.
The longest perimeter small octagon was computed in \cite{Audet2007a},
but its perimeter was already stated to four correct digits 
in \cite{Griffiths1975}.

The extreme polygons for other powers of $2$ are still unknown.

\section{Long perimeter small triacontadigons}

We now consider the case $n=32$. It is known that a convex small $32$-gon
with maximal perimeter exists and possesses a full set of $32$ diameter chords
of unit length.
The diameter graph consists of an odd cycle of unit segments supplemented 
by hanging diameters attached to the vertices of the cycle. 
The following table summarizes the
longest perimeters of small $32$-gons available so far.
For comparison, the upper bound is
\begin{align*}
    2\cdot 32 \sin\frac\pi{2\cdot 32} =3.14033115|\mcg{6954753}.
\end{align*}

\begin{table}[h]
	\centering
	\begin{tabular}{cD{.}{.}{1.16}c} \toprule 
	cycle&\multicolumn{1}{c}{perimeter} & source \\ \midrule
	3&3.1403|\mcg{23421103532}&\cite{Tamvakis1987} \\
	17&3.140331|\mcg{085836778}& \cite{Mossinghoff2006}, $B_{32}^*$ \\
	23&3.14033115|\mcg{4141625}& \cite{Bingane2021}, $D_{32}^*$ \\
	21&3.14033115|\mcg{6355381}& \cite{Xue2021} \\
	21&3.14033115|\mcg{6954614}& new\\
	\bottomrule
	\end{tabular}
	\caption{Long perimeters of small $32$-gons.} 
	\label{fig:r32}
\end{table}

Based on our computational approach, which will be described elsewhere, we 
claim with reasonable certainty
that the last row actually contains the \notion{longest perimeter of a small
triacontadigon}. We also have enough evidence to happily confirm that the 
small hexadecagon $D_{16}^*$
of \cite{Bingane2021} realizes the \notion{extreme} perimeter $3.13654771648661$
for $n=16$.

\section{Appendix}

\begin{figure}[htb]
    \begin{subfigure}[b]{.46\linewidth}
        \centering\includegraphics*[width=73mm]{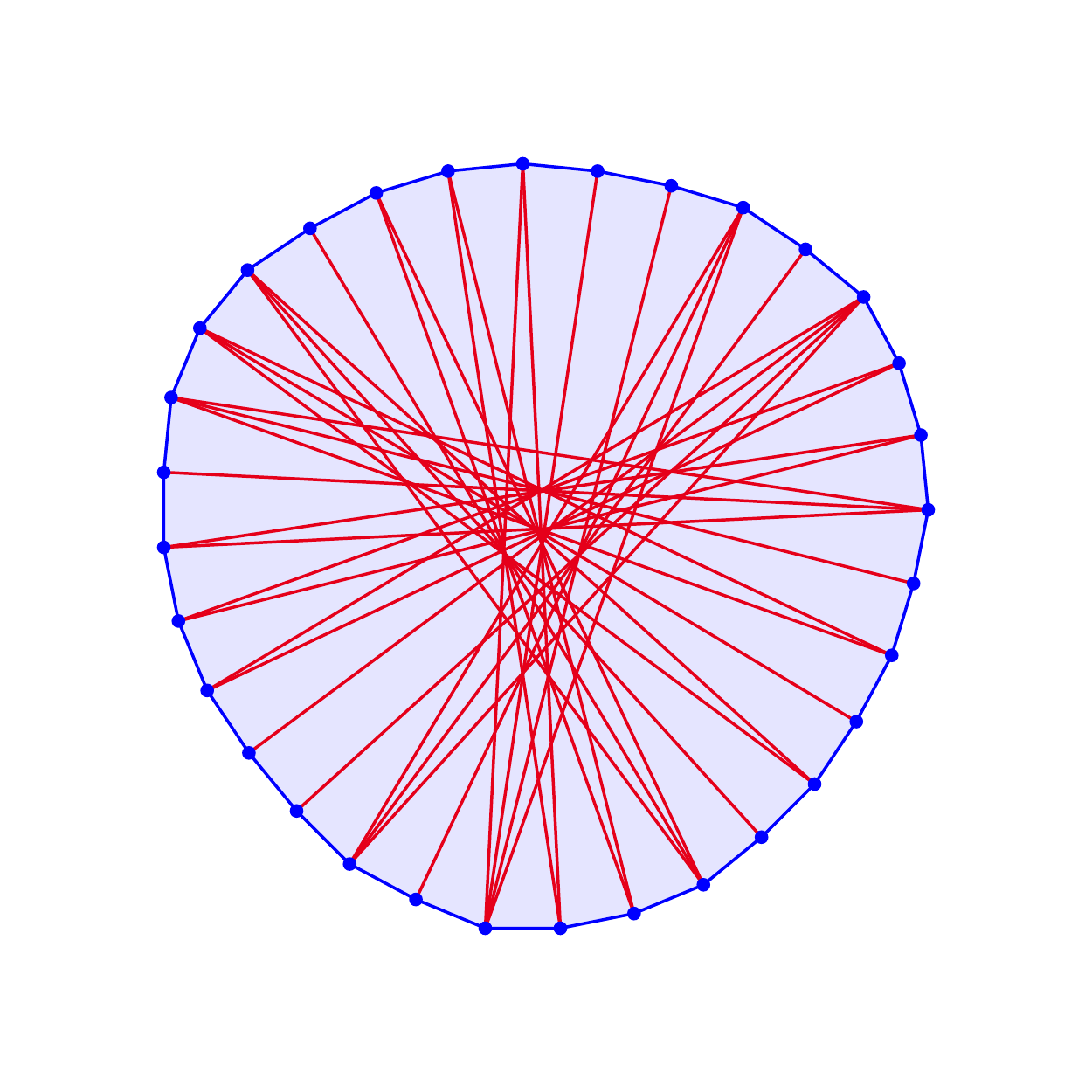}
    \end{subfigure}
    \begin{subfigure}[b]{.46\linewidth}
        \centering\includegraphics*[width=73mm]{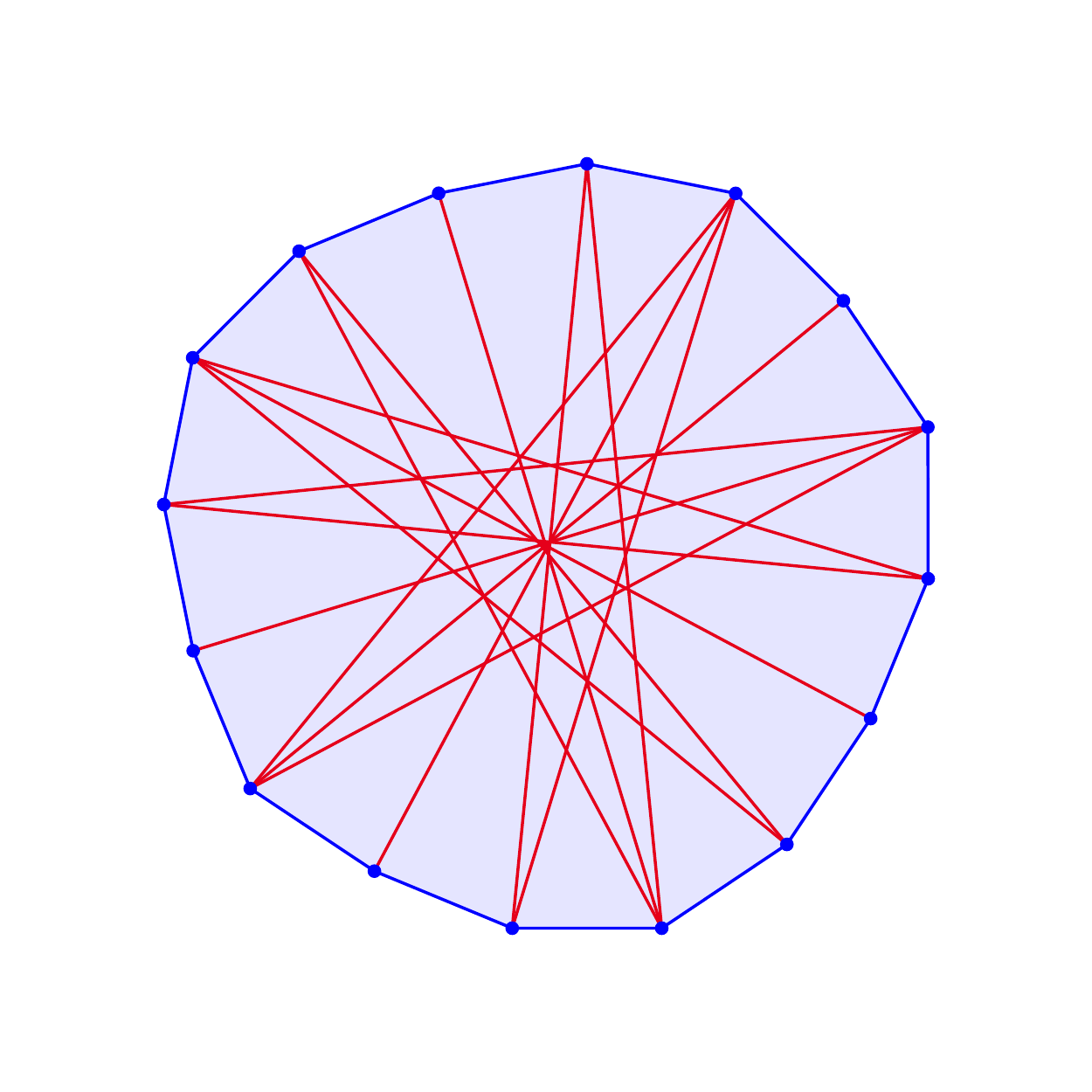}
    \end{subfigure}
    \caption{The longest perimeter small triacontadigon
	and hexadecagon.} 
    \label{fig:opt}
\end{figure}
The plot shows the longest perimeter small $32$-gon
and $16$-gon including their diameter graphs.
To allow for the verification of the stated properties we provide the cartesian
coordinates of the vertices of our $32$-gon in the following table. 

\begin{table}[h]
	\centering
	\begin{tabular}{D{.}{.}{1.16}D{.}{.}{1.16}}
 0.098134910054388 & 0. \\
 0.1943840604354086 & 0.01914514251146182 \\
 0.28504868084122 & 0.056699649798818 \\
 0.3609076612541907 & 0.1189555734331473 \\
 0.4302992695838125 & 0.1883471408623677 \\
 0.4848199067577719 & 0.2699429769467566 \\
 0.5310802950708347 & 0.3564898244762278 \\
 0.5595673636200257 & 0.4503986731664625 \\
 0.5787126617485553 & 0.5466475302839486 \\
 0.569093927380149 & 0.6443096742296724 \\
 0.5406070476551426 & 0.7382188490444155 \\
 0.4943467071289093 & 0.824766185605239 \\
 0.418487098558063 & 0.887022659054982 \\
 0.3368905726741924 & 0.941543807819504 \\
 0.2429806029380585 & 0.970031147229744 \\
 0.1467305815533126 & 0.989176494078293 \\
 0.04906745502719413 & 0.998795466978675 \\
 -0.04859580696564962 & 0.989176473983678 \\
 -0.142506038633965 & 0.960689119789753 \\
 -0.2290542902439535 & 0.914428123805245 \\
 -0.3106513808547955 & 0.859906619698457 \\
 -0.3729082997819818 & 0.7840463942725833 \\
 -0.4104634795827763 & 0.6933804898426428 \\
 -0.4200827044600564 & 0.595717036499123 \\
 -0.4200828511732027 & 0.4975810104628068 \\
 -0.4009376289104606 & 0.4013307043569394 \\
 -0.3633828226693777 & 0.310664976005411 \\
 -0.3088616948467671 & 0.2290680548412778 \\
 -0.2466051794150348 & 0.1532080706018036 \\
 -0.1772128446516963 & 0.0838156012797867 \\
 -0.0906652207709911 & 0.03755484577878503 \\
 0. & 0. \\
\end{tabular}
\caption{Coordinates of the vertices of the triacontadigon.} 
\label{fig:gon32}
\end{table}

\printbibliography

\end{document}